\documentclass[12pt]{article}
\usepackage[amsmath]{e-jc}

\usepackage{graphicx}

\usepackage{caption}
\usepackage{wrapfig}
\usepackage{enumerate}
\usepackage{fge}

\usepackage[numbers]{natbib}

\dateline{Jan 24, 2022}{Mar 28, 2023}{Apr 21, 2023}


\Copyright{The authors. Released under the CC BY-ND license (International 4.0).}

\title{On enumeration and entropy of ribbon tilings}

\author{Yinsong Chen \quad Vladislav Kargin \\
\small Department of Mathematics and Statistics\\[-0.8ex]
\small Binghamton University\\[-0.8ex] 
\small Binghamton, U.S.A.\\
\small\tt \{ychen276, vkargin\}@binghamton.edu
}

\newcommand{\bal}[1]{\begin{align*}#1\end{align*}}

\newcommand{\ovln}[1]{\overline{#1}}

\newcommand{\R}{\mathbb{R}}
\newcommand{\E}{\mathbb{E}}
\newcommand{\Z}{\mathbb{Z}} 
\newcommand{\A}{\mathcal A}

\newcommand{\TT}{\mathcal{T}} 
\newcommand{\eps}{\varepsilon}

\newcommand{\sm}{\fgebackslash} 
\newcommand{\Ext}{\mathrm{Ext}}

\theoremstyle{plain}
\newtheorem{theo}{Theorem}[section]

\theoremstyle{definition}
\newtheorem{defi}[theo]{Definition}

\begin{document}

\maketitle


\begin{abstract}
The paper considers ribbon tilings of large regions and their per-tile entropy (the logarithm of the number of tilings divided by the number of tiles). For tilings of general regions by ribbon tiles of length $n$, we give an upper bound on the per-tile entropy as $n - 1$. For growing rectangular regions,  we prove the existence of the asymptotic per-tile entropy and show that it is bounded from below by $\log_2 (n/e)$ and from above by $\log_2(en)$. For growing generalized ``Aztec Diamond'' regions and for growing ``stair'' regions, the asymptotic per-tile entropy is calculated exactly as $1/2$ and $\log_2(n + 1) - 1$, respectively.
\end{abstract}

\section{Introduction}
 Let a region $R \in \R^2 $ be a union of finite number of unit squares $[k, k+1]\times [l, l+1]$, with $k, l \in \Z$.  Assume that the interior of $R$ is connected and simply connected. We consider tilings of such regions by ribbon tiles. 
 
%

\begin{defi} 
\label{defi_ribbon_tile}
A \emph{ribbon} of length $n$, or an \emph{$n$-ribbon}, is a connected sequence of $n$ unit squares, each of which (except the first one) comes directly above or to the right of its predecessor.
\end{defi}

\begin{figure}[htbp]
\begin{minipage}[b]{0.45\linewidth}
\centering
              \includegraphics[height = 0.7\linewidth]{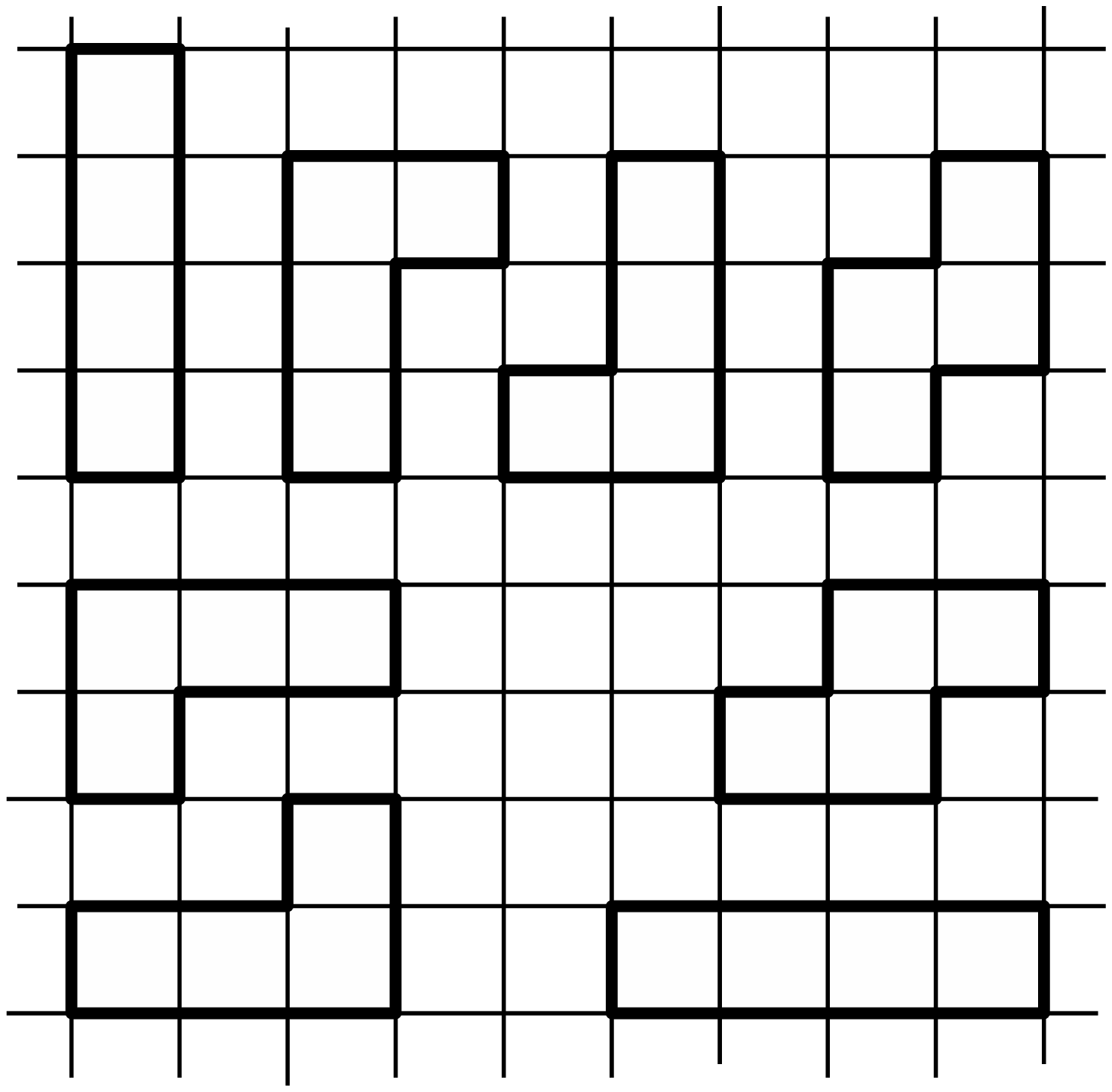}
              \caption{Eight types of $4$-ribbons.}
              \label{figRibbons}
\end{minipage}
\hspace{0.5cm}
\begin{minipage}[b]{0.45\linewidth}
\centering
              \includegraphics[width= 0.8\linewidth]{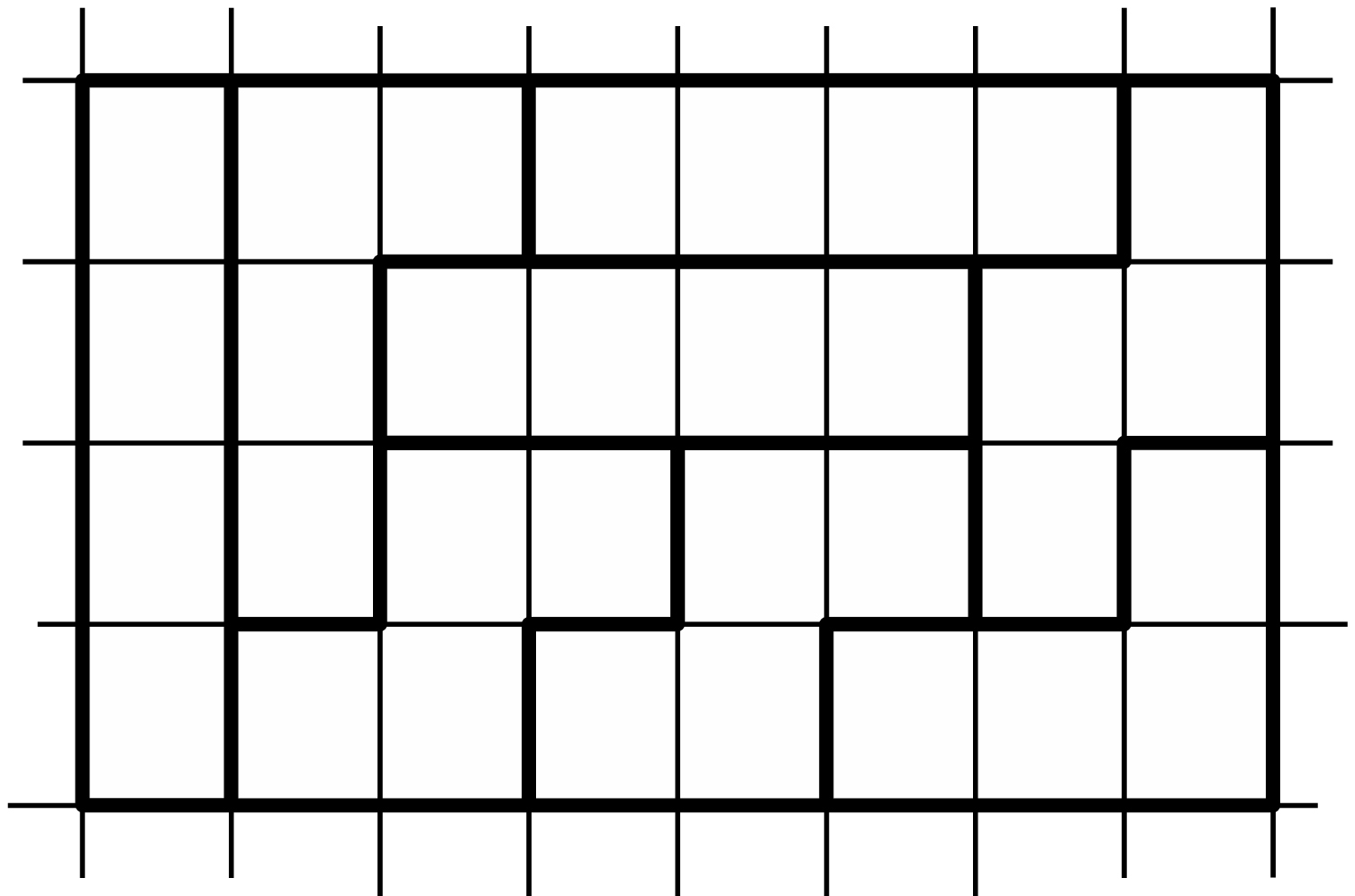}
              \caption{A ribbon tiling of a $4 \times 8$ rectangle.}
              \label{figTiling4by8}
\end{minipage}
\end{figure}

These objects are also called \emph{border strips} or \emph{rim hooks} in the literature.\footnote{See \cite{stanley99}, Section 7.17, p.345} See an illustration in Figures \ref{figRibbons} and \ref{figTiling4by8}.
 Dominoes are a particular case of ribbon tiles with $n = 2$. Ribbon tiles with $n = 3$ are called (right-oriented) $180$-trominoes in \cite{agmsv2020}.\footnote{However, in \cite{agmsv2020} strictly horizontal or strictly vertical tiles are excluded, so there are four types of 3-ribbons but only two types of right-oriented $180$-trominoes.} The study of ribbon tilings for $n \geq 3$ was initiated in  \cite{pak2000a}, and developed extensively in \cite{sheffield2002a}.

 Typical questions about tilings are:
 \begin{enumerate}
 \item Is it possible to tile a region $R$? 
 \item How many different tilings exist?
 \item What is the distribution of tile shapes in a typical tiling?
 \item How to sample random tilings efficiently?
 \end{enumerate}
 
The existence question was studied  in \cite{sheffield2002a}, who gave an algorithm for checking if a simply-connected region is tileable by $n$-ribbons.   The algorithm is linear in the area of the region. In \cite{agmsv2020}, it was shown that for general regions (which are allowed to be non-simply connected with arbitrary number of holes), the existence of tilings by $180$-trominoes  is an $NP$-complete decision problem.

In this paper, we focus on question (2), the question of enumeration. For enumeration problems, much is known about domino tilings and lozenge tilings in the triangular lattice,  and precious little is known about any other types of tilings.  An exception is tilings by $T$-tetrominoes, where the enumeration of tilings has been related to the evaluation of the Tutte graph polynomial at the argument pair $(3, 3)$ (see \cite{korn_pak2004} and \cite{merino2008}). In addition, some numerical results have been obtained in \cite{hutchinson_widom2015} for octagonal tilings. The goal of this paper is to rectify to a certain extent this deplorable scarcity of enumeration results. 

Let us define the \emph{per-tile entropy} of $n$-ribbon tilings of a region $R$ as the binary logarithm of the number of tilings divided by the 
number of ribbons in each tiling, that is, 
\bal{
\mathrm{Ent}_n(R) = \frac{\log_2(|\mathcal{T}_n(R)|)}{T},
}
where $\mathcal{T}_n$ is the set of all $n$-ribbon tilings of the region $R$, and $T$ is the number of tiles in each tiling (that is, $T=Area(R)/n$).

Suppose we consider a sequence of regions $\big(R_t\big)_{t=1}^\infty$  with 
$A_t = \mathrm{area}(R_t)\to  \infty$ as $t \to \infty$. Then we are interested in the existence and the value of the limit 
\bal{ 
\mu^{(n)}(R_t) =  \lim_{t \to \infty}\mathrm{Ent}_n(R_t) 
=\lim_{t \to \infty} \frac{\log_2\big| \TT_n(R_t)\big|}{A_t/n}.
}

For domino tilings of rectangles and of Aztec diamonds, this limit can be calculated explicitly due to formulas obtained in \cite{Kasteleyn1961} and \cite{Temperley_Fisher1961} in the case of rectangles, and in \cite{eklp1992c}  in the case of Aztec diamonds.  
In particular, for a fixed $N$, and $M \to \infty$, the formulas for the number of domino tilings imply that 
 \bal{
\lim_{M \to \infty} \frac{\log_2\big| \TT(R_{N,M})\big|}{(NM)/2} = \frac{2}{N}\sum_{l = 1}^{\lfloor N/2 \rfloor}
  \log_2\bigg(\cos \frac{l \pi}{N + 1} + \sqrt{1 + \cos^2 \frac{l \pi}{N + 1}}\bigg)
 }

If both $N$ and $M$ approach $\infty$, then the limit is 
\bal{
\mu^{(2)}(``Rectangle") = \frac{2 G} { \pi \ln 2} = 0.841\,266\,940\, 7 \ldots,
} 
where $G$ denotes the Catalan constant ($G = 1 - 3^{-2} + 5^{-2} - 7^{-2} + \ldots$).

For the tilings of the ``Aztec Diamond'', we have 
   \bal{
          \mu^{(2)}(``AD") = \frac{1}{2} < \mu^{(2)}(``Rectangle"). 
    }

More generally, let $R$ be a region in $\R^2$ bounded by a piecewise smooth, simple closed curve $\partial R$, without cusps.
Approximate $R$ by a region $R_\eps = R \cap \eps \Z^2$. 
Then, it was shown in \cite{cohn_kenyon_propp2001a} that the asymptotic growth in the number of domino tilings of  $R_\eps$, $\eps \to 0$, and therefore the entropy of the set of tilings, 
can be described using a certain functional on height functions associated with domino tilings.

For ribbon tilings, it is not difficult to calculate that the number of tilings of an $n\times n$ square by $n$-ribbons is $n!$ (see Lemma \ref{lemmaNTilingsRectangle} below). If we let $n$ grow then we obtain the entropy of $\log_2(n!)/n \sim \log_2(n) -\log_2(e)$. However, in this limit we let both the region size and the ribbon length  grow. It would be more natural to have the length of the ribbons fixed and the size of the region growing to infinity. 

A glimpse of what might happen in this case can be obtained from a result in  \cite{alexandersson_linus2019}, which gives the number of $n$-ribbon tilings for $n \times 2n$ rectangle. Let us denote this number $a_n := \big|\TT_n(R_{n, 2n})\big|$. The formula is 
\bal{
a_n = \frac{1}{2}\sum_{i = 1}^n \frac{i(n - i + 1)}{n+2} \binom{n - 1}{i - 1} \binom{n +3}{i +1} a_{i - 1} a_{n - i}, 
}
with the initial condition $a_0 = 1$. For example, $a_1 = 1$, $a_2 = 5$, $a_3 = 61$, $a_4 = 1379$. (This is sequence A115047 in OEIS.) 

Somewhat mysteriously, these numbers coincide with Weil-Peterson volumes of moduli spaces of algebraic curves, and for their asymptotic expression we have:
\bal{
a_n = \frac{(2n)!}{C^n}\Big(c_1 - \frac{c_2}{n} + \ldots \big), 
}
where $C = 2.496918339\ldots$ is a constant that can be expressed in terms of Bessel functions and their derivatives (see formula (0.9) in \cite{kmz1996}). 

So, if $n \to \infty$, we have 
\bal{
\frac{\log_2 a_n}{2n} &\sim \log_2 (2n) - \log_{2} e - \frac{1}{2} \log_2 C 
\\
&= \log_2 n - \log_2 e + 1 - \frac{1}{2}\log_2 C,
}
which shows that changing the $n\times n$ square to the $n \times 2n$ rectangle leads to a significant increase in the limiting per tile entropy by  $1 - \frac{1}{2}\log_2 C = 0.339926\ldots $. However, in this example we still have the situation in which both the region size and the ribbon length grow. In our considerations below, we will focus on the limit in which the ribbon length is fixed while the region size grows to infinity.

 Other enumeration results for ribbon tilings were also obtained in \cite{stanley99} and \cite{stanley2002}. These results hold for ribbon tilings that are allowed to have ribbons of varying length, which is different from the situation we consider here. We explain Stanley's results and compare them with our results in Appendix. 
 
 Finally, a recent paper \cite{richter2023} considers \emph{coverings} of rectangles by monotonous polyominoes which are close relatives of our ribbon tiles.   In coverings, as distinct from tilings, the tiles can overlap and the focus of \cite{richter2023} is on the minimal number of tiles needed to cover a rectangle rather than on the number of tilings. 
 
 In this paper, we prove that for an arbitrary region the per-tile entropy for $n$-ribbon tilings is bounded from above by $n - 1$. For growing rectangular regions, we prove that the per-tile entropy converges to a finite limit and we bound the limit from below by $\log_2(n/e)$ and from above by $\log_2(en)$. This property, that the per-tile entropy increases in the length of the ribbon tile is not universal, since we show that for tilings of generalized Aztec diamonds, the entropy is always $1/2$, for all values of $n$. 

Then, in Theorem \ref{theoStairEntropy} we calculate the exact value of the asymptotic entropy for certain stair regions provided that  ribbon tiles have odd length $n$. In this case, the asymptotic entropy equals $\log_2(n + 1) - 1 = \log_2(n) - 1 + o(1)$. In another study (which is going to be published separately), we show that in the case of thin rectangles of height $M = n$ and growing width $N$, the asymptotic entropy is bounded from above by $\log_2 n$. Together these results suggest that it would be interesting to calculate 
\bal{
\lim_{n \to \infty} \big[\mu^{(n)}(R_t) - \log_2(n)\big]
}
exactly for a family of rectangles $R_t$ that grows both in width and length.

The rest of the paper is organized as follows. In Section \ref{sectionResults}, we explain our results in more detail.  In Section \ref{sectionProofs}, we provide proofs. And in Appendix, we describe Stanley's enumeration results for ribbon tilings  and compare them with our findings. 

%

\section{Results}
\label{sectionResults}

First, we establish a simple general upper bound on the number of $n$-ribbon tilings. 
\begin{theo} 
\label{theoUpperBound}
Let $R$ be an arbitrary simply connected region of the square lattice which is a union of $nT$ squares. Then the number of tilings of $R$ by $n$-ribbons is bounded from above by $2^{(n - 1)T}$. In particular, for arbitrary sequence of regions $\big(R_t\big)_{t = 1}^\infty$ with growing area $A_t$,
\bal{
\limsup_{t \to \infty} \frac{\log_2 \big|\TT_n(R_t)\big|}{A_t/n} \leq n-1
}
\end{theo}

For rectangular regions, the existence result for the entropy was shown by Yinsong Chen in his Ph.D. thesis (\cite{chen2020}). We provide a proof for the reader's convenience.  

Let $R_{M, N}$ denote a rectangle with $M$ rows and $N$ columns and let $\big|\mathcal{T}_{M, N}\big|$ be the number of $n$-ribbon tilings of this region. 

\begin{theo}
\label{theoExistence}
For a sequence of  $(R_t)_{t = 1}^\infty$ of $M_t \times N_t$ rectangles, assume that both $M_t$ and $N_t \to \infty$ as $t \to \infty$ and that rectangles $R_{M_t, N_t}$ are all tileable.  Then
the following limit exists and finite:
\bal{
 \mu^{(n)}(R_t) := \lim_{t\to \infty} \frac{\log_2 \big|\mathcal{T}_{M_t, N_t}\big|}{T_t} 
 = \sup_{M, N} \frac{\log_2\big|\mathcal{T}_{M, N}\big|}{T} 
\leq n - 1,
}
where 
 $T = MN/n$  denote the number of ribbons in each tiling of the rectangle $R_{M,N}$. 
\end{theo}
%
%
%

We can further give a lower bound on the entropy. 
\begin{theo}
\label{theoLowerBound}
For a sequence of  $(R_t)_{t = 1}^\infty$ of $M_t \times N_t$ rectangles, assume that both $M_t$ and $N_t \to \infty$ as $t \to \infty$ and that rectangles $R_{M_t, N_t}$ are all tileable.   Then,
\bal{
\mu^{(n)}(R_t) \geq
 \log_2 n  - \log_2 e 
+ n^{-1}\Big(
\frac{1}{2} \log_2 n + \log_2 \sqrt{2\pi}
\Big).
}
\end{theo}

Observe that while the lower bound is growing in $n$, this bound and the upper bound in Theorem \ref{theoUpperBound} are far apart. However, for rectangular regions we can give a better upper bound than that in Theorem \ref{theoUpperBound}.

\begin{theo}
\label{theoUpperBoundRect}
For a sequence of  $(R_t)_{t = 1}^\infty$ of $M_t \times N_t$ rectangles, assume that both $M_t$ and $N_t \to \infty$ as $t \to \infty$ and that rectangles $R_{M_t, N_t}$ are all tileable. Then,
\bal{
\mu^{(n)}(R_t) \leq
 \log_2 n  +\log_2 e .
}
\end{theo}

Besides these results, we also have exact expressions for the tiling entropy of some non-rectangular regions. 


 \begin{figure}[htbp]
\begin{minipage}[b]{0.45\linewidth}
\centering
              \includegraphics[height=\linewidth]{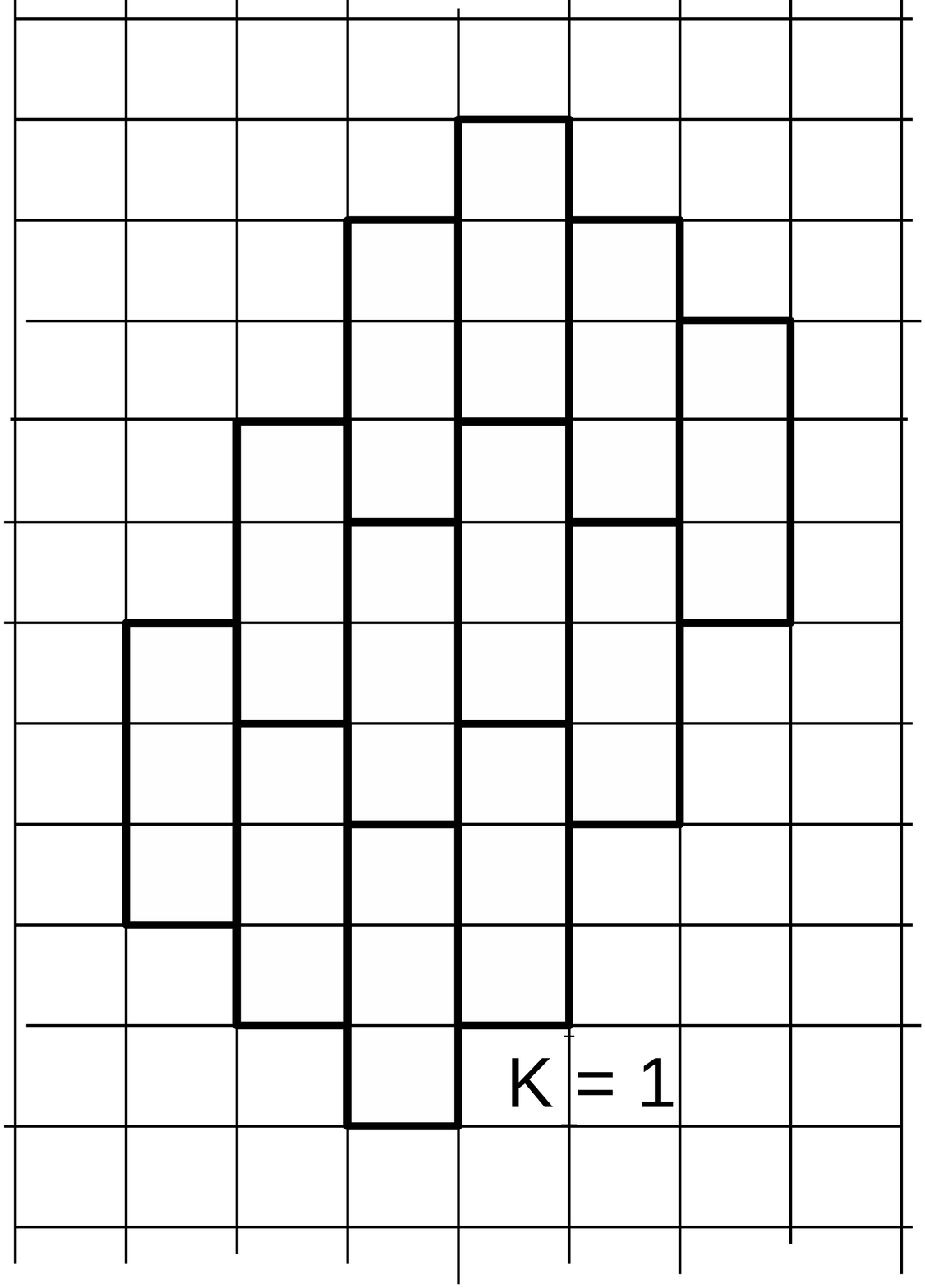}
              \caption{A tiling of $AD(3, 3, 1)$ by 3-ribbons.}
              \label{figGeneralizedAztec}
\end{minipage}
\hspace{0.1cm}
\begin{minipage}[b]{0.45\linewidth}
\centering
              \includegraphics[height=\linewidth]{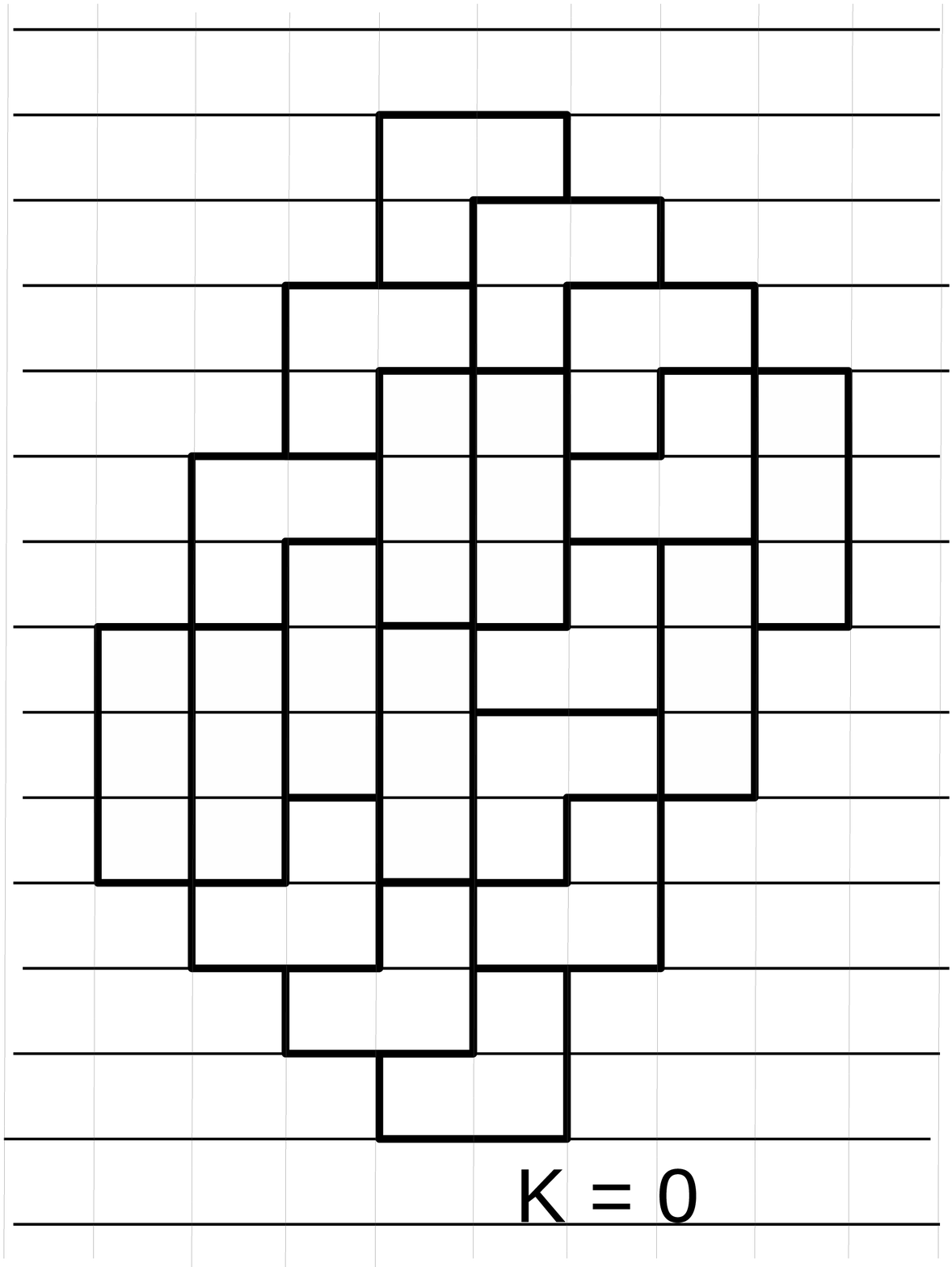}
              \caption{A tiling of $AD(4, 3, 0)$ by 3-ribbons.}
              \label{figGenAztecTiling}
\end{minipage}
\end{figure}

First, note that usual Aztec diamonds are non-tileable by $n$-ribbons if $n \geq 3$. (For the case $n = 3$, see Theorem 17 in \cite{agmsv2020}.)  For this reason, we consider generalized Aztec diamonds which are tileable by $n$-ribbons. We define a generalized Aztec diamond region $AD(N, n, k)$ as shown in Figure \ref{figGeneralizedAztec}. The parameter $N$ measures the size of the region, -- the number of columns equals $2N$, and the two longest columns can be covered by $N$ ribbon tiles.  The parameters $k$ and $n$ satisfy inequalities $0 \leq k \leq n - 2$, and $n$ corresponds to the length of ribbons which will be used to tile the region, while $k$ is the ``offset'' in the diamond shape.  That is, $k$ measures the amount by which one of the largest columns is shifted relative to the other one. The region $AD(N, 2, 0)$ is the usual Aztec diamond. 
 
 
  Note that  $AD(N, n, k)$ can be tiled by $N(N + 1)$ of $n$-ribbons. Two ribbon tilings of generalized Aztec diamonds are shown for illustration in the Figures \ref{figGeneralizedAztec} and \ref{figGenAztecTiling}.
   
 \begin{theo}
 \label{theoAztec}
 The number of tilings of region $AD(N, n, k)$ by $n$-ribbons equals  
 $2^{N(N+1)/2}$, and as $N \to \infty$, the limit per-tile entropy for this sequence of these regions is $\mu^{(n)} = 1/2$. 
 \end{theo} 

 In particular, the per-tile entropy is not zero but it does not depend on the length of ribbon tiles, in contrast with results for the rectangular regions. 
 
 \textbf{Remark.} The proof of this Theorem in Section \ref{sectionEntropyAD} essentially builds a bijection between domino tilings of the Aztec Diamond $AD(N, 2, 0)$ and $n$-ribbon tilings of $AD(N, n, k)$. Examples suggest that under this bijection vertical dominoes are mapped to vertical $n$-ribbons and horizontal dominoes are mapped to ribbon tiles which can have two possible types. Both of these types are vertical except for exactly one horizontal step. In particular, one can observe in random tilings of  $AD(N, n, k)$  an analogue of the Aztec circle effect characteristic for domino tilings of the regular Aztec Diamond.

 \begin{figure}[htbp]
\begin{minipage}[b]{0.45\linewidth}
\centering
              \includegraphics[width=\linewidth]{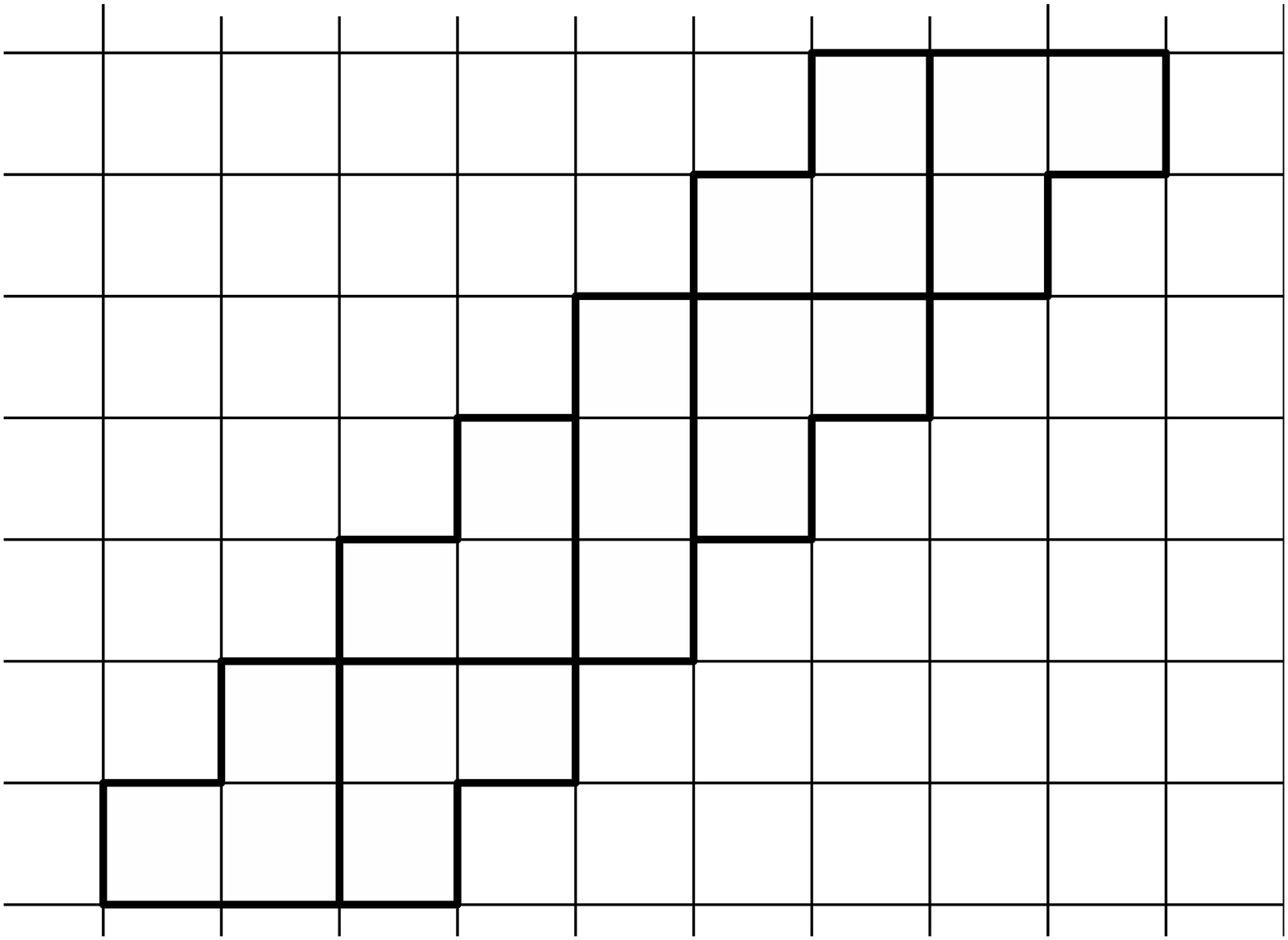}
              \caption{A tiling of a stair $St_7^{(3)}$by $3$-ribbons}
              \label{figStair7by3}
\end{minipage}
\hspace{0.5cm}
\begin{minipage}[b]{0.45\linewidth}
\centering
              \includegraphics[width=\linewidth]{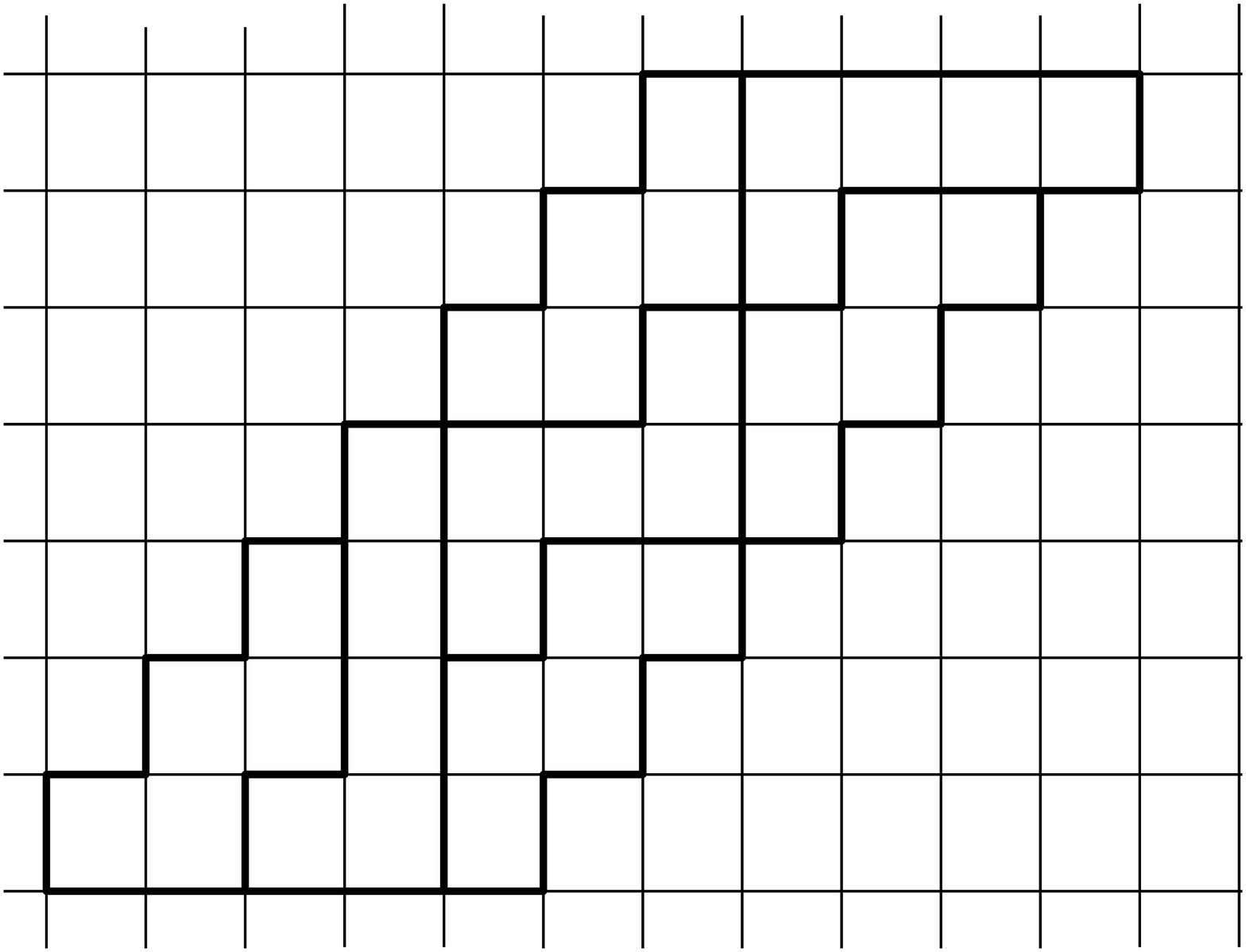}
              \caption{A tiling of a stair $St_7^{(5)}$by $5$-ribbons}
              \label{figStair7by5}
\end{minipage}
\end{figure}

Finally, let us define the stair region $St_M^{(n)}$ of size $M \times n$ as shown in Figures  \ref{figStair7by3} and \ref{figStair7by5}.     That is, a stair $St_M^{(n)}$ has $M$ rows, the length of each row is $n$, and each row is offset by 1 square to the right relative to the row below. 

\begin{theo}
\label{theoStairEntropy}
Let $\{St_M^{(n)}\}\Big|_{M = 1}^\infty$ be a sequence of  $M\times n$ stairs. For every odd $n$, the limit tiling entropy for this family of regions equals $\log_2(n + 1) - 1$.
\end{theo}

In this case the entropy $\mu^{(n)}$ is growing with $n$ as $\log_2(n)$. This is similar to what is observed in the case of rectangular regions.


%
\section{Proofs}
\label{sectionProofs}

\subsection{Preliminaries}
First of all, let us introduce some notation. 
Let $s_{xy}$ denote the square whose south-west corner has coordinates $(x, y)$. We will say that a \emph{level} of a square $s_{xy}$ is $l(s_{xy}) =  x + y$. The \emph{root square} in a tile is the square with the smallest level.  The \emph{level of the tile} is the level of its root square. 
 
 It was proved in \cite{sheffield2002a} that, for a given region $R$, each tiling has the same number of tiles in a specific level. In particular, in each level we can enumerate tiles
 from left to right (from the tile with the smallest $x$-coordinate of its root square to the tile with the largest $x$-coordinate). Then, let $t_{l, i}$, $i = 1, \ldots, k_l$ denote the $i$-th tile in the level $l$ in this enumeration. Here $k_l$ depends only on $R$. This enumeration gives us an unambiguous way to refer to a specific tile in any tiling of region $R$. 
 
 Next, let us describe a couple of constructions from \cite{sheffield2002a}. First, every simply connected region $R$ can be put in correspondence with a graph $G_R$ whose vertices are identified with tiles $t_{l,i}$ and, additionally, with squares on the border of the region. We will describe below how the edges are defined. Some edges in the graph $G_R$ are endowed with an orientation which depends only on region $R$ but not on the tiling. These edges are called forced. This construction gives a graph $G_R$ with a partial orientation $\tau_R$. A second construction shows that each ribbon tiling of $R$ determines a complete acyclic orientation $\tau$ on the edges of this graph which is in agreement with $\tau_R$.
 
 One of the central results in \cite{sheffield2002a}, is that every
acyclic orientation $\tau$ on the graph $G_R$ that extends the partial orientation $\tau_R$  can be realized by a ribbon tiling of region $R$. This gives a bijection between acyclic orientations on $G_R$ extending $\tau_R$ and ribbon tilings of $R$.
 

 For details of these constructions, see \cite{sheffield2002a}. Here we want to briefly explain to the reader how to build graph $G_R$ and partial orientation $\tau_R$, and how the tilings can be translated to orientations on $G_R$. We somewhat simplify the definitions keeping in mind the examples in which we are interested in. 
 
 Recall that tile $t_{l, i}$ is the $i$-th tile in the level $l$.
 We use $t_{l,i}$ also as labels for vertices in the graph $G_R$.  In order to handle the border conditions, we assume that there is a fixed tiling of the region outside of $R$ and include in the graph $G_R$ the vertices corresponding to tiles outside but adjacent to the border of $R$.\footnote{As it turns out, in most of our examples, these vertices can be omitted without any change in the number of admissible acyclic orientations. It is possible to introduce other border conditions and study, for example, the tilings of a torus instead of a rectangle. However, we will not discuss such extensions in this paper.}  the border conditions do not matter to us.) We postulate that there is an edge between two different $t_{l,i}$ and
 $t_{l', i'}$ if and only if $|l - l'| \leq n$. This postulate defines graph $G_R$ for a  region $R$. 

  A  tiling of $R$ defines an orientation on each edge of this graph in the following way. Consider two ribbon tiles $t$ and $t'$ in a tiling, and imagine that every square in these tiles projects light in the north-west direction. 
 Then if some light from $t$ is absorbed by $t'$ than we say that $t'$  is \emph{to the left} of $t$, and orient the 
 corresponding edge in the graph $G_R$ from $t'$ to $t$.  As a result we obtain an orientation of the graph $G_R$ and it is clear that this orientation is acyclic.
 
 Consider the orientation of edges between vertices in the same level. An edge $(t_{l,i}, t_{l,i'})$ is oriented from $t_{l,i}$ to $t_{l,i'}$ if and only if $i < i'$. These orientations are obviously the same for every tiling and we include them in the partial orientation $\tau_R$ of the graph $G_R$.
 
In addition, it was proved in \cite{sheffield2002a} that if  $l - l' = \pm n$, then the orientation 
of an edge  $(t_{l, i}, t_{l',i'})$ is always the same for every ribbon tiling of region $R$. This orientation is completely determined by the geometry of the region. Therefore, we also include these orientations in the partial orientation $\tau_R$, and we say that the edges with orientations in $\tau_R$ are \emph{forced edges}. It is easy to build $\tau_R$ by looking at one specific tiling of region $R$. 

The question of  enumerating ribbon tilings of $R$ is reduced, therefore, to the questions of enumerating acyclic orientations of $G_R$, which are in agreement with partial orientation $\tau_R$.

%

\subsection{An upper bound on the entropy}
\label{sectionUpperBound}
\begin{proof}[Proof of Theorem \ref{theoUpperBound}] 
 
 There are $2^{n-1}$ different shapes for an $n$-ribbon, since each type can be encoded by a sequence of  $0$s and $1$s that has $n - 1$ elements. In this sequence $0$ corresponds to the ribbon continuing to the left and $1$ to the continuing up. For example, sequences $0 0 \ldots 0$ and $1 1 \ldots 1$ encode the horizontal and vertical $n$-ribbons, respectively.  
 
Each ribbon tiling of a region $R$ can be mapped to an assignment which says what is the  shape of a tile $t_{l, i}$. 

It is clear that given such an assignment, the tiling can be unambiguously recovered. Indeed, proceed level by level in the order of increasing $i$. If all previous tiles have been already placed in the region, we can determine where the root square of the current tile $t_{l,i}$ is located, and then the shape of this tile, read from the assignment, determines the placement of the entire tile in the region. 

Therefore, the number of tilings is no greater than the number of possible shape assignments which is $(2^{n - 1})^T$.
\end{proof}

%

\subsection{The existence of the entropy limit for rectangular regions}

\label{sectionLimitExistence}

\begin{proof}[Proof of Theorem \ref{theoExistence}]
Let $(M, N) \in \Omega$ mean that rectangle $R_{M, N}$ is tileable by ribbon tiles of length $n$. It is easy to check that this holds if and only if at least one of $M$ and $N$ is divisible by $n$. It is also easy to check that rectangles $R_{M, N}$ and $R_{N, M}$ have the same number of tilings. 
Let
\bal{
S = \sup_{(M, N)\in \Omega} \frac{\log_2 \big|\TT_{M, N}\big|}{MN},
}
For any $\eps > 0$, choose $(\ovln{M},\ovln{N}) \in \Omega$,  such that $\log_2 \big|\TT_{\ovln M, \ovln N}\big|/(\ovln M \ovln N) > S - \eps$. Without loss of generality we can choose the rectangle so that $n|\ovln M$.

We are going to show that for all $M$ and $N$ sufficiently large, $\log_2 \big| \TT_{M, N}\big|/( M  N) > S - \eps$. By our observation about transposed rectangles above, it is enough to prove this claim for the case when $n | M$.

So,   let $M = p \ovln M + k$ and $N = q \ovln N + l$, where $1 \leq k \leq \ovln M$ and $1 \leq l \leq \ovln N$. Note that $k$ is divisible by $n$.

Then the rectangle $R_{M, N}$ can be split in $pq$ rectangles congruent to $R_{\ovln M, \ovln N}$, and 3 rectangles $R_{p\ovln M, l}$, $R_{k, q \ovln N}$ and $R_{k, l}$. All of these rectangles are tileable and by super-additivity of the number of tilings we get:
\bal{
\frac{\log_2 \big|\TT_{M, N}\big|}{MN} \geq 
\frac{p q\log_2 \big|\TT_{\ovln M, \ovln N}\big|}{(p\ovln M + k)(q  \ovln N + l)} 
= \frac{\log_2 \big|\TT_{\ovln M, \ovln N}\big|}{(\ovln M + k/p)(\ovln N + l/q)}\geq S - \eps 
}
for all sufficiently large $p$ and $q$. 

Hence the limit over increasing $M$ and $N$ exists and equals the supremum. By Theorem \ref{theoUpperBound} this limit is finite and $ \leq n - 1$. 
\end{proof}


\subsection{A lower bound on the entropy for rectangular regions}

\label{sectionLowerBoundRectangular} 

\begin{lemma}
\label{lemmaNTilingsRectangle}
The number of tilings of an $n\times N$ rectangle by $n$-ribbons is equal to $N!$ for $N \leq n$ and to $\frac{1}{2}N!$ for $N = n + 1$.  
\end{lemma}
\begin{proof}
For $N \leq n$, every tiling of an $n \times N$ rectangle by $n$-ribbons  has one tile in each of the levels $0, \ldots n - 1$. The Sheffield graph $G_R$ of this region is the complete graph on $N$ vertices that correspond to the tiles of the tiling and all its edges are free. (For rectangular regions, the vertices corresponding to border squares do not impose additional restrictions for acyclic orientations and can be safely ignored.) Hence the number of tilings of this region equals the number of acyclic orientations on the complete graph $K_N$, which equals $N!$ (the number of vertex orderings).

If $N = n + 1$, then the graph is again a complete graph on $N$ vertices but the orientation on one edge is forced  (the vertices in levels $0$ and $n$, $t_0$ and $t_n$, are comparable but $t_0$ must be to the left of $t_n$). That means that exactly half of all $N!$ possible orientations of the graph have the correct orientation on this edge. It follows that there are $(n + 1)!/2$ tilings of the $n \times (n + 1)$ rectangle by $n$-ribbons. 
\end{proof}

\begin{lemma}
We have the following lower bound for the limit per-tile entropy of tilings of an $n \times N$ rectangles by $n$-ribbons: 
\bal{
\mu^{(n)}_n \geq \frac{\log_2(n!)}{n} 
\geq \log_2 n  - \log_2 e 
+ n^{-1}\Big(
\frac{1}{2} \log_2 n + \log_2 \sqrt{2\pi}
\Big),
}
\end{lemma}
\begin{proof}
The claim  follows from Lemma \ref{lemmaNTilingsRectangle} and the super-additivity for the logarithm of the number of tilings.  For a large $N$ we split the $n\times N$ strip in $t = \lfloor N/n \rfloor$ of $n\times n$ squares and a remainder region. By the previous lemma and super-additivity, the number of tilings is greater than $(n!)^t$. Hence,  the logarithm of the number of tilings divided by the number of tiles $N$, 
\bal{
\frac{\log_2 \big|\mathcal{T}_{n, N}\big|}{N} \geq \frac{\log_2(n!) \lfloor N/n \rfloor}{N}
}
By taking the limit $N \to \infty$ we demonstrate the first inequality in the corollary. 
  The second inequality follows from a well-known lower estimate on $\log(n!)$, see, for example, formula (1.53) on p. 17 in \cite{spencer2014}.
\end{proof}

\begin{proof}[Proof of Theorem \ref{theoLowerBound}]
As in the proof of the previous lemma, the claim  of the theorem follows from the super-additivity of the logarithm of the number of tilings. In this case we divide the rectangle in $\lfloor M/n \rfloor$ strips of the size $n \times N$, plus a remaining strip, and use the estimate in the proof of the previous lemma to bound the logarithm of the number of tilings of each strip from below by 
$
 \lfloor N/n \rfloor \log_2(n!).
$
Then the logarithm of the total number of tilings is bounded from below by 
\bal{
\lfloor N/n \rfloor \times \lfloor M/n \rfloor \log_2(n!)
}
After dividing by the number of tiles $N M/ n$ and taking the limit, we find the claimed inequality for the entropy.
\end{proof}


\subsection{An upper bound on the entropy for rectangular regions}

\label{sectionUpperBoundRectangular} 
If $H$ is a subgraph of $G_R$, let $\A(H)$ be the set of all acyclic orientations on $H$ which agree with the partial orientation $\tau_R$. Also for shortness, we will write \emph{a.o.} for acyclic orientations that agree with the partial orientation $\tau_R$.

  \begin{figure}[htbp]
\begin{minipage}[b]{0.45\linewidth}
\centering
              \includegraphics[width=\linewidth]{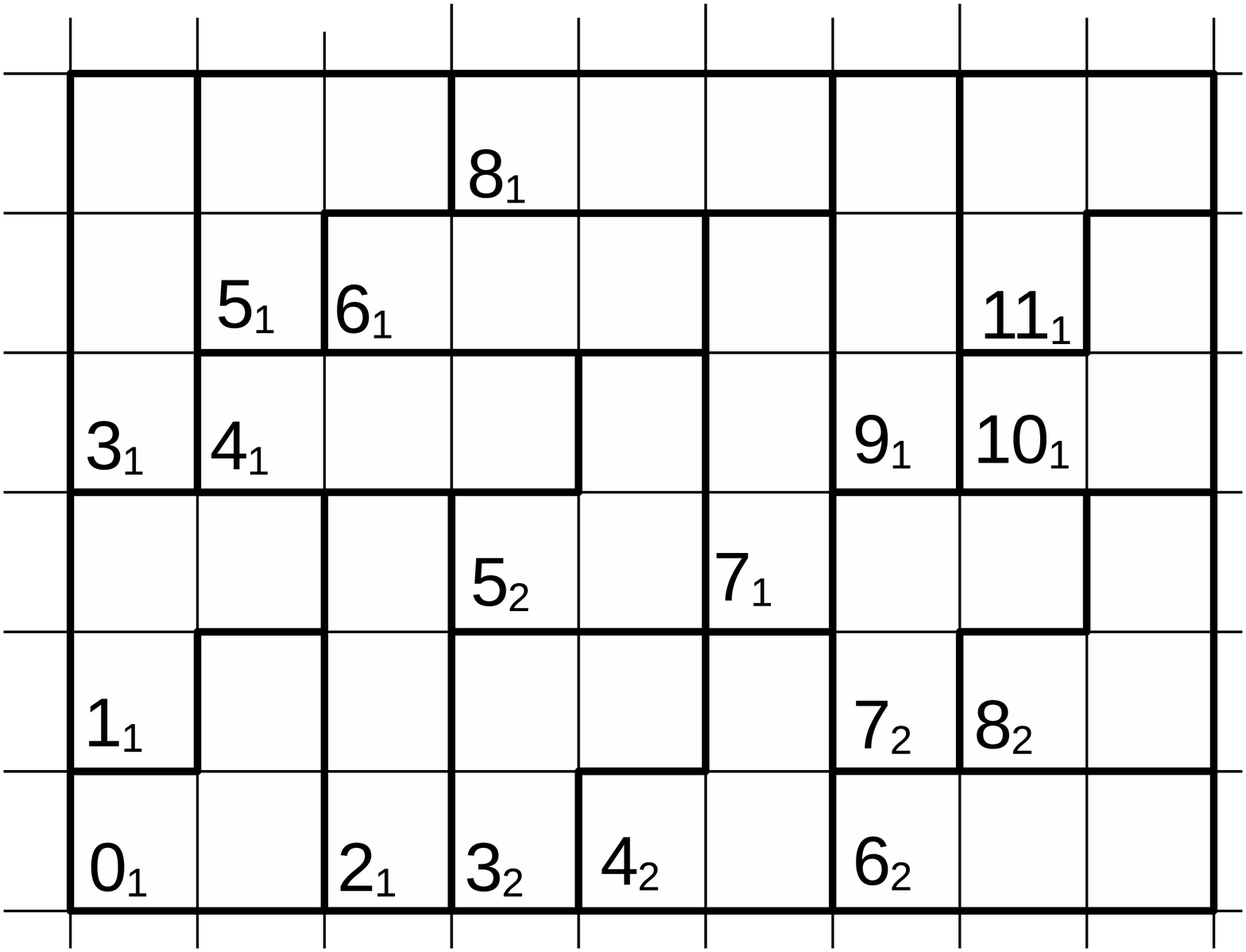}
              \caption{A tiling of a $6\times 9$ rectangle by $3$-ribbons, with ribbons labeled by the level. The subscript is the order number of the ribbon in the level.}
              \label{figRectangle}
\end{minipage}
\hspace{0.5cm}
\begin{minipage}[b]{0.45\linewidth}
\centering
              \includegraphics[height=0.8\linewidth]{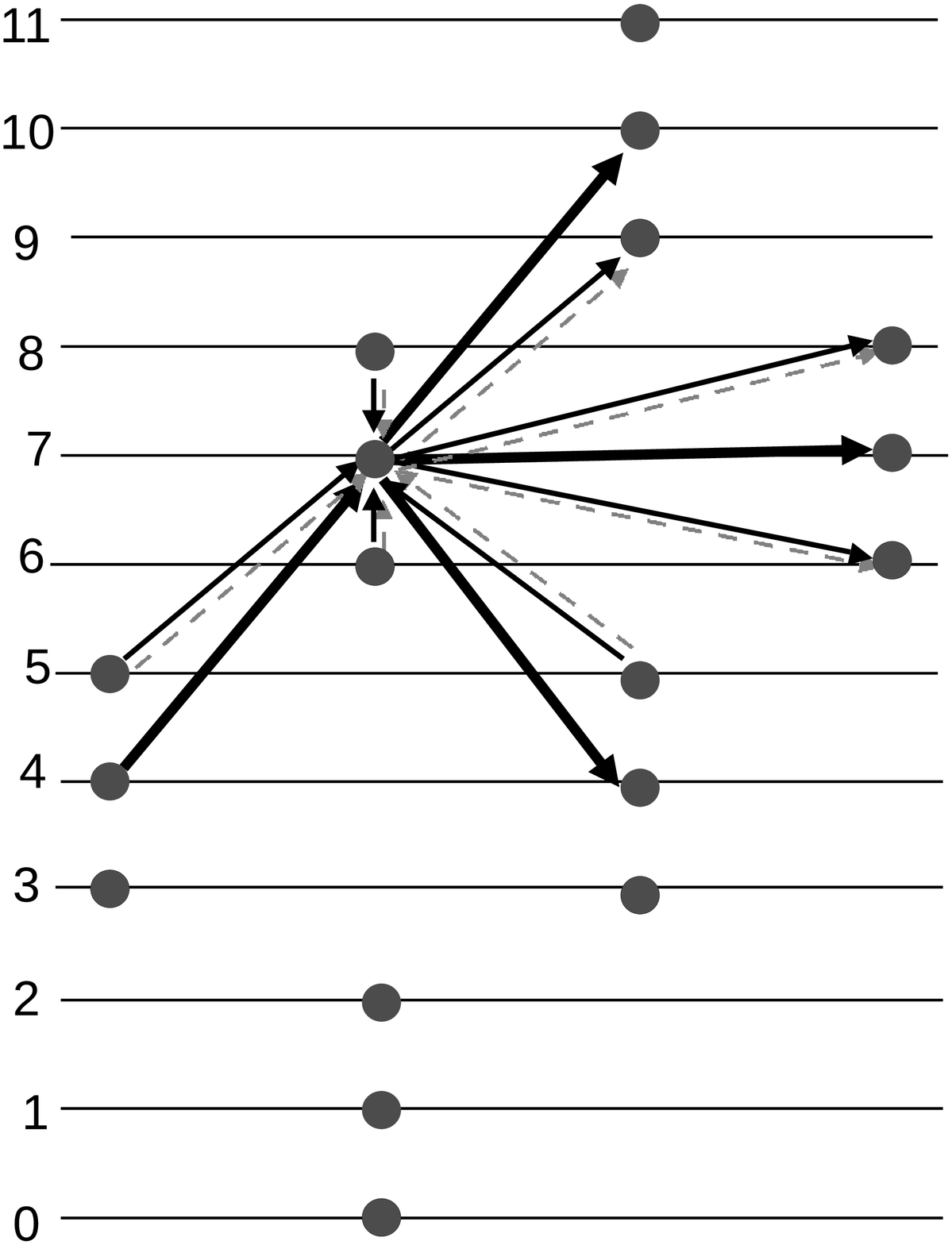}
              \caption{The edges of the Sheffield graph of the $6\times 9$ rectangle incident to the vertex $7_1$. The forced edges are shown in solid, the other edges are doubled by a dashed line, with the orientation induced by the tiling in Figure \ref{figRectangle}.}
              \label{figRectangle_graph}
\end{minipage}
\end{figure}

Consider an $M \times N$ rectangle $R$. It can be shown that it is tileable by $n$-ribbons if and only if either $M$ or $N$ is divisible by $n$. Assume without loss of generality that $n | M$. 
An illustration of a 3-ribbon tiling for a rectangle  is shown in Figure \ref{figRectangle}.

Let
 $H_l$ denote the subgraph of $G_R$ induced by vertices $v$ with $l(v) \leq l$, where $0 \leq l \leq  l_{max}$ and $l_{max} = M + N - n - 1$ is the maximum possible level of a tile in an $n$-ribbon tiling of $R$. (For example, $l_{max} = 11$ in Figures \ref{figRectangle} and \ref{figRectangle_graph}.)
 
Define growth factor
\bal{
g_l := \frac{|\A(H_l)|}{|\A(H_{l - 1})|} 
}
 for $ 1 \leq l \leq l_{max}$, with the convention that $|\A(H_{0})| = 1$, so that $g_1 = |\A(H_1)| = 2$.  
 
 Let $T_l := |\{v: l(v) = l\}$, the number of vertices at level $l$, and $S_l := \sum_{k = l - n + 1}^l T_l$, the number of vertices at levels between $l - n + 1$ and $l$, inclusive.
 
 \begin{lemma}
 \label{lemmaGrowthRate}
 For all $1 \leq l \leq l_{max}$, 
 \bal{
 g_l \leq \binom{S_l}{T_l}.
 }
 \end{lemma}
 Before proof, let us introduce some additional notation. If $\beta$ is an a.o. on a graph $G$ and $H \subset G$ is a subgraph of $G$, then we write $\beta|_H$ to denote the \emph{restriction} of $\beta$ to $H$. Obviously, $\beta|_H$ is an a.o. on $H$. 
 
 If $\alpha \in \A(H)$ and $H \subset G$, then we call orientation $\beta$ on $G$ an \emph{extension} of $\alpha$ if $\beta \in \A(G)$ and $\alpha = \beta|_H$.  Clearly, an extension is determined by $\alpha$ and the orientations on edges in $E(G) \sm E(H)$, where $E(G)$ and $E(H)$ denote the sets of edges of $G$ and $H$, respectively.  These orientations must be such that no directed cycle is created in $\beta$. We denote the set of all extensions of $\alpha \in \A(H)$ from $H$ to $G$ as $\Ext(\alpha; H, G)$.
 
 \begin{proof}[Proof of Lemma  \ref{lemmaGrowthRate}]
 By restriction, every a.o. $\beta$ in $\A(H_l)$ corresponds to an a.o. $\alpha$ in $\A(H_{l - 1})$, $\alpha = \beta|_{H_{l - 1}}$. To prove the lemma, it is enough to show that for every $\alpha \in \A(H_{l - 1})$ there are no more than $\binom{S_l}{T_l}$ extensions of $\alpha$ in $\Ext(\alpha, H_{l - 1}, H_{l})$. 
 
 Let 
 \bal{
 V_l &= \{v \in G_R: l(v) = l\}, \\
 U_l &= \{v \in G_R: l(v) \in [l - n + 1, l)\}, \\
 \ovln U_l & = U_l \cup V_l. 
 }
 By abusing notation, we will also use $V_l$, $U_l$, $\ovln U_l$ to denote the subgraphs of $G_R$ induced by respective sets of vertices. 
 
 Now, if $\alpha \in \A(H_{l - 1})$ and $\beta \in \Ext(\alpha, H_{l - 1}, H_l)$, then 
 \bal{
 \beta|_{\ovln U_l} \in \Ext(\alpha|_{U_l}, U_l , \ovln U_l).
 }
 Moreover, 
 \bal{
 E(H_l) \sm E(H_{l - 1}) = \big(E(\ovln U_l) \sm E(U_{l})\big) \cup F, 
 }
 where $F$ is the set of edges between vertices in $V_l$ and vertices in $V_{l - n}$. Since the orientations on $F$ are forced by partial orientation $\tau_R$, the orientations assigned by extension $\beta|_{\ovln U_l}$  of $\alpha|_{U_l}$ to edges in $\E(\ovln U_l) \sm \E(U_{l})$ also completely determine the extension $\beta$ of $\alpha$. It follows that 
 \bal{
 |\Ext(\alpha, H_{l - 1}, H_l)| \leq |\Ext(\alpha|_{U_l}, U_l , \ovln U_l)|.
 }
 In particular, in order to prove the lemma, it is enough to show that $|\Ext(\alpha, U_l , \ovln U_l)|\leq \binom{S_l}{T_l}$ for every $\alpha \in \A(U_l)$.
 
 Next, note that $U_l$ and $\ovln U_l$ are complete graphs. (Indeed, the levels of any two vertices in $\ovln U_l$ differ by no more than $n - 1$, hence they are comparable with respect to the left-of relation, and therefore connected by an edge.) Also, note that the only edges with forced orientation in $\ovln U_l$ are the edges between vertices in $V_l$. 
 
 On a complete graph, acyclic orientations are in one-to-one correspondence with linear orders on vertices. Now, an a.o. $\alpha$ on $U_l$ determines a linear order on vertices of $U_l$ and the forced orientations on $V_l$ determine a linear order on  vertices of $V_l$. Then, the number of linear orders on $\ovln U_l$ consistent with given linear orders on $U_l$ and $V_l$ equals the number of ways to insert ordered vertices of $V_l$ between ordered vertices of $U_l$. By an elementary combinatorial formula, this number equals 
 \bal{ 
 \binom{|U_l| + |V_l|}{|V_l|} = \binom{S_l}{T_l}. 
 }
 We showed that for every $\alpha \in \A(U_l)$, $|\Ext(\alpha, U_l , \ovln U_l)|= \binom{S_l}{T_l}$ and by observations above, this completes the proof of the lemma.
 \end{proof}
 
 Let $T_{max} = \max \{T_l, 0\leq l \leq l_{max}\}$ and 
 \begin{equation}
 \label{def_L}
 L = \max \{l: T_l = T_{max}\}.
\end{equation}
  In other words, $L$ is the highest level among those that have the most vertices. For example, in Figure \ref{figRectangle_graph}, $T_{max} = 2$ and $L = 8$.
 
 \begin{lemma}
 \label{lemma_growth_estimate}
 If $l \leq L$, then $g_l \leq (e n)^{T_l}$.  
 \end{lemma}
 Here $e$ is the base of the natural logarithm. 
 \begin{proof}
 By a well-known property of binomial coefficients (see, for example, inequality (5.14) on p. 59 in \cite{spencer2014}),  
 \bal{
 \binom{S_l}{T_l} \leq \Big(e \frac{S_l}{T_l}\Big)^{T_l}. 
 }
 Observe that $(T_l)$ is a non-decreasing sequence for $0 \leq l \leq L$, and therefore $S_l \leq n T_l$ in this range. Then, the conclusion of the lemma follows from Lemma  \ref{lemmaGrowthRate}.
 \end{proof}

 
\begin{proof}[Proof of Theorem \ref{theoUpperBoundRect}]
Let $H$ and $H'$ be subgraphs of $G_R$ induced by vertices with $l(v) \leq L$ and vertices with $l(v) \geq L - n$, respectively, where $L$ is as defined in (\ref{def_L}). Then, it is clear that $E(G_R) \subset E(H) \cup E(H')$. 

For $\alpha \in \A(G_R)$, we can define the map $\alpha \mapsto (\alpha|_H, \alpha|_{H'})$. This is an injective map from  $\A(G_R)$ to $\A(H) \times \A(H')$ (since $\alpha$ can be recovered unambiguously from $\alpha|_H$ and $\alpha|_{H'}$), and therefore 
\bal{
|\A(G_R)| \leq |\A(H)| |\A(H')|. 
}
Then, by using Lemma  \ref{lemma_growth_estimate}, we have
\bal{
|\A(H)| = \prod_{l = 1}^L g_l \leq \prod_{l = 1}^L (en)^{T_l} = (en)^{\sum_{l = 1}^L T_l} = (en)^{|H|}.
}
Similarly, by using the fact that graph $G_R$ is symmetric, we can obtain the estimate $|\A(H')| \leq (en)^{|H'|}$, and therefore,
\bal{
|\A(G_R)| \leq (en)^{|H| + |H'|}. 
}
Then,
\bal{
|H| + |H'|&= |G_R| + |H \cap H'|\leq |G_R| + (n + 1) T_{max} \\
 &\leq |G_R| + (n + 1) \frac{M}{n} \leq |G_R| + 2 M.
}
Note that $|G_R| = MN/n$. Then, we have:
\bal{
\mu^{(n)}(R_t) &= \lim_{M, N \to \infty} \frac{\log_2(|\A(G_R)|)}{MN/n}
\\
& \leq 
\lim_{M, N \to \infty} \frac{\log_2(en) (MN/n + 2 M)}{MN/n}
= \log_2(en). 
}

\end{proof}

%
\subsection{Number of tilings of a generalized Aztec Diamond} 

\label{sectionEntropyAD}
 
 \begin{proof}[Proof of Theorem \ref{theoAztec}]
It is easy to check that the Sheffield graph $G_R$ and the partial orientation $\tau_R$ of the generalized Aztec diamond $AD(N, n, k)$ for $n$-ribbon tilings are isomorphic to the Sheffield graph and the partial orientation of $AD(N, 2, 0)$ for $2$-ribbon tilings, that is,  with that of domino tilings of the standard Aztec diamond. This implies that these graphs have the same number of acyclic orientations that agree with the partial orientation. Hence, the number of $n$-ribbon tilings of $AD(N, n, k)$ equals  the number of domino tilings of the standard Aztec diamond, and this number was computed in the celebrated result in \cite{eklp1992c}.
 \end{proof}
 
  Note that this proof provides a bijection between domino tilings of $AD(N, 2, 0)$ and $n$-ribbon tilings of $AD(N, n, k)$. In this bijection, two tilings correspond to each other if they induce the same acyclic orientation on the isomorphic Sheffield graphs of $AD(N, 2, 0)$ and $AD(N, n, k)$. Intuitively, one can think about this bijection as that one can judiciously add $n - 2$ squares to each domino tile of the domino tiling so as to lengthen the shape $AD(N, 2, 0)$ vertically and make it coincide with $AD(N, n, k)$.
 

%

\subsection{Exact value for the limit entropy of a stair region} 
\label{sectionEntropyStair}

The proof is based on the exact enumeration of the number of ribbon tilings for these regions. 

\begin{theo}
\label{theoStair}
The number of tilings of the stair region $St_M^{(n)}$ by $n$-ribbons is given by the following formulas. If $n$ is odd, then 
\bal{
|\TT_n(St_M^{(n)})| = 
\begin{cases}
M! & \text{ for } M \leq \frac{n + 1}{2}, 
\\
\Gamma\big(\frac{n + 1}{2}\big)  \big(\frac{n + 1}{2}\big)^{M - (n-1)/2},
& \text{ for }  M > \frac{n + 1}{2}.
\end{cases}
}
If $n$ is even, then $|\TT_n(St_M^{(n)})|$ equals the number of tiling of an $n/2 \times M$ rectangle by $n/2$-tiles, namely  
\bal{
|\TT_n(St_M^{(n)})|= |\TT_{n/2}(R_{n/2, M})|.
}
\end{theo}


 \begin{figure}[htbp]
\begin{minipage}[b]{0.45\linewidth}
\centering
              \includegraphics[width=\linewidth]{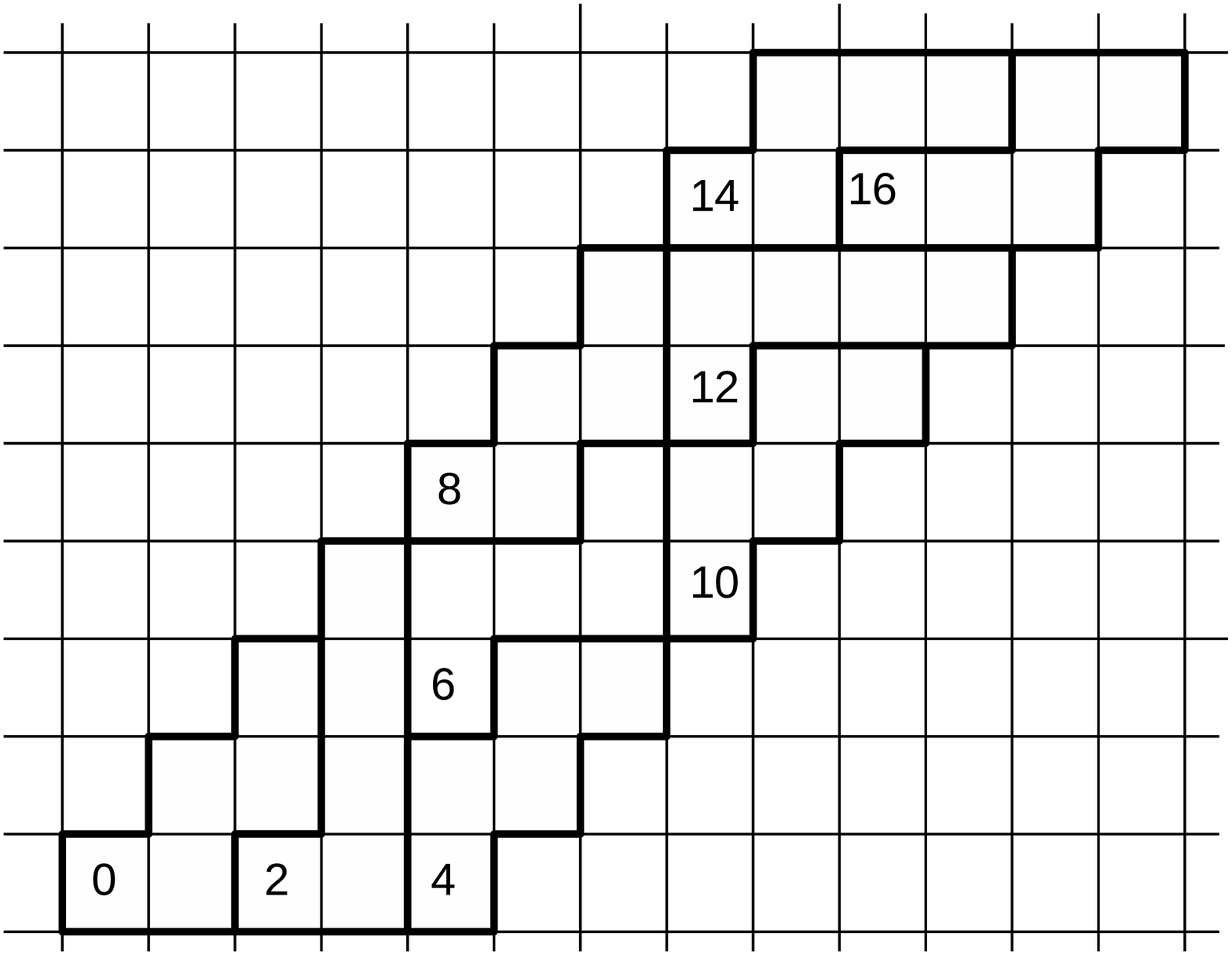}
              \caption{A tiling of $St_9^{(5)}$ by $5$-ribbons, with ribbons labeled by the ribbon level.}
              \label{figSt9_5}
\end{minipage}
\hspace{0.5cm}
\begin{minipage}[b]{0.45\linewidth}
\centering
              \includegraphics[width=\linewidth]{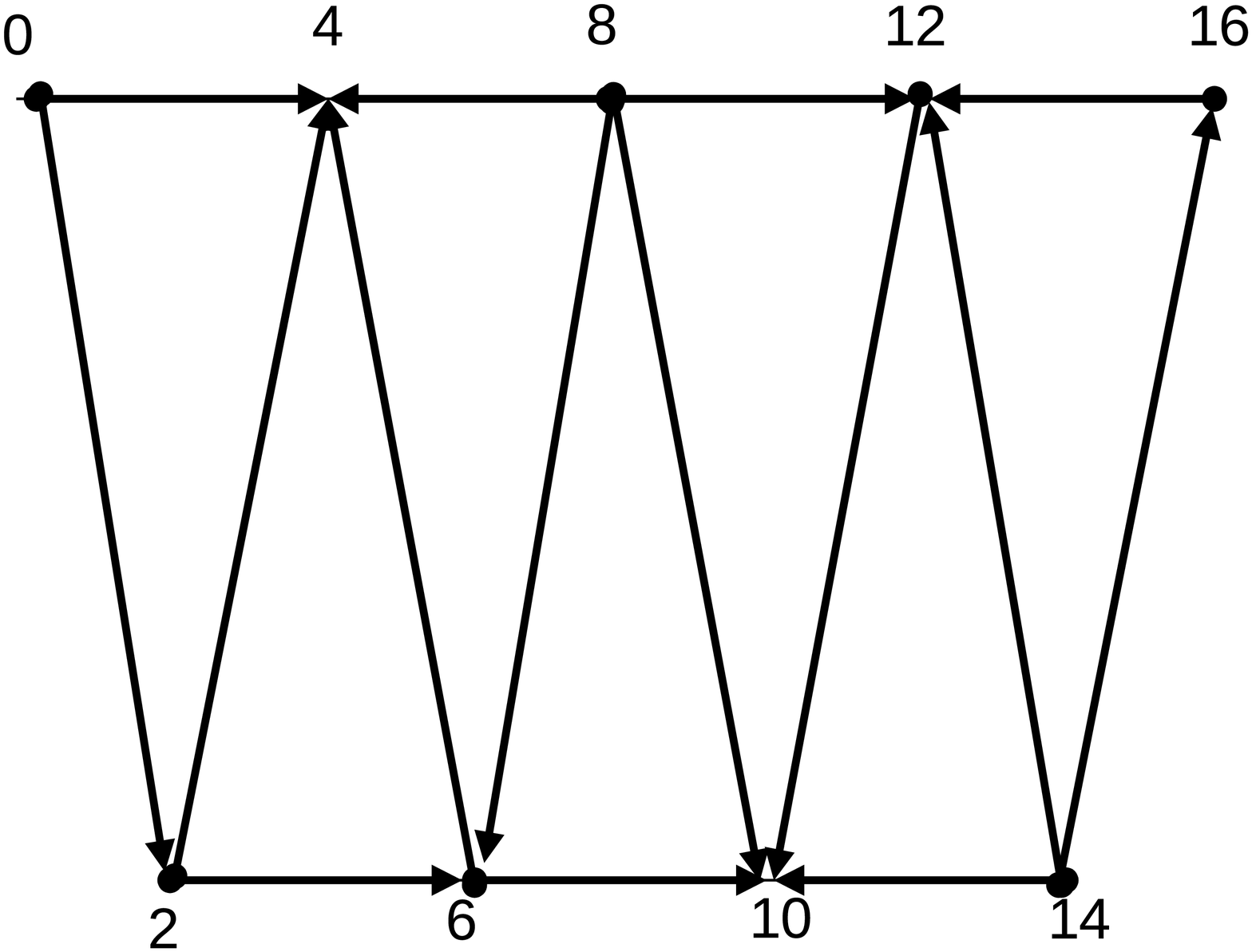}
              \caption{The Sheffield graph of $St_9^{(5)}$ with the orientation induced by the tiling in Figure \ref{figSt9_5}.}
              \label{figSt9_5_graph}
\end{minipage}
\end{figure}

\begin{proof}

Consider the case when $n$ is odd. Every tiling of the region $St_M^{(n)}$ by $n$-ribbons has $M$ tiles, one tile in each level $0$, $2$, $\ldots$,  $2(M - 1)$. We denote these tiles $t_0, t_2, \ldots , t_{2(M - 1)}$. Then, the Sheffield graph $G_{St_M^{(n)})}$ that corresponds to this region  has edges in the following list, 
\bal{
\big(t_{2k}, t_{2(k + i)}\big) , \text{ where } k = 0, \ldots, M - 2 \text{ and } i = 0, \ldots, \frac{n - 1}{2}, 
}
provided that both end-points are well defined.  (See an example for $n = 5$ in Figures \ref{figSt9_5} and \ref{figSt9_5_graph}). 

Crucially, this graph does not have any edges with forced orientation, since (i)~the differences between levels of tiles are even and therefore are different from $n$, and (ii)~ every level has no more than one tile. 

By Sheffield's theorem, the number of tilings of the region $St_M^{(n)}$ is equal to the total number of acyclic orientations of the graph $G_{St_M^{(n)}}$, since the partial orientation of $G_{St_M^{(n)}}$ is empty. If $M \leq (n + 1)/2$, then the graph  $G_{St_M^{(n)}}$ is the complete graph $K_M$ and the number of acyclic orientations is $M!$. 

If $M > (n + 1)/2$, then we use  Stanley's theorem (\cite{stanley73}) that the number of acyclic orientations of a graph $G$ equals to $|\chi_G(-1)|$, where $\chi_G(t)$ denotes the chromatic polynomial of the graph $G$. The chromatic polynomial for the graph $G_{St_M^{(n)})}$ is calculated in Lemma \ref{lemmaChromatic}, and we get the following formula for the number of acyclic orientations:
\bal{
\Big| (-1) ( - 2) \ldots  (- m ) \big(-( m + 1)\big)^{M - m}\Big|,
}
where $m = (n - 1)/2$. This proves the theorem for the odd $n$.

For the case when $n$ is even, it is easy to check that the graph of the $St_M^{(n)}$ for $n$-ribbon tilings coincide with the graph of  $R_{n/2, M}$ for $n/2$-ribbon tilings. This implies that the number of tilings is the same. 
\end{proof}

\begin{lemma}
\label{lemmaChromatic}
Let $n$ be odd, let $m = (n - 1)/2$, and let graph $\Gamma_M = G_{St_M^{(n)}}$ be as defined in the proof of Theorem \ref{theoStair}, with $M > (n + 1)/2 = m + 1$. Then, the chromatic polynomial of the graph $\Gamma_M$ is 
\begin{equation}
\label{equChromatic}
\chi_{\Gamma_M}(\lambda) = \lambda (\lambda - 1) \cdots (\lambda - m + 1) (\lambda - m)^{M - m}.  
\end{equation}
\end{lemma}

\begin{proof}
We use induction on $M$ and Read's theorem from \cite{read1968} that says that if a graph $G$ is a union of two subgraphs $X$ and $Y$, which overlap in a complete graph $K_s$ on $s$ nodes, then the chromatic polynomial of the graph $G$  is 
\begin{equation}
\label{equRead}
\chi_G(\lambda) = \frac{\chi_X(\lambda) \chi_Y(\lambda)}{\chi_{K_{s}}(\lambda) } = \frac{\chi_X(\lambda) \chi_Y(\lambda)}{\lambda^{(s)}}, 
\end{equation}
where $\lambda^{(s)}$ is the factorial monomial:
\bal{
\lambda^{(s)} = \lambda (\lambda - 1) \cdots (\lambda - s + 1)
}
In our case the graph $X$ is the restriction of graph $\Gamma_M$ to vertices $t_0$, 
$t_2$, \ldots, $t_{2(m + 1)}$, and the graph $Y$ is the restriction of graph $\Gamma_M$ to vertices $t_2$, \ldots, $t_{2(M - 1)}$.  They intersect in complete graph $K_{m + 1}$.
(In the example in Figure \ref{figSt9_5_graph}, $m = 2$,  $X = \Gamma_M | \{ 0, 2, 4, 6\} $ and  $Y = \Gamma_M | \{ 2, 4, 6, 8, 10, 12, 14, 16\}$.)


The graph $X$ equals the complete graph on $m + 2$ vertices with the edge $e = (t_0, t_{2(m + 1)})$ removed. A well known property of chromatic polynomials relates the polynomial of a graph to the polynomials of the graphs obtained by a contraction and a removal of an edge, respectively. Namely, for every simple graph $G$ and for all $e \in E(G)$, 
\bal{
\chi_G(\lambda) = \chi_{G\sm e}(\lambda) - \chi_{G/ e}(\lambda),
}
where $G\sm e$ denotes $G$ with the edge $e$ deleted, and $G/ e$ denotes $G$ with the edge $e$ contracted to a point. 

By applying this property to the complete graph $K_{m + 2}$, note that the removal of an edge $e$ leads to the graph $X$ and a contraction of $e$ leads to $K_{m + 1}$. Hence, 
\bal{
\chi_X(\lambda) &= \chi_{K_{m + 2}}(\lambda) +
\chi_{K_{m + 1}}(\lambda) 
\\
&= \lambda \cdots (\lambda - m) (\lambda - m - 1) +
\lambda \cdots (\lambda - m) 
\\
&=  \lambda \cdots (\lambda - m + 1) (\lambda - m)^2. 
}
For the base of the induction, we note that if $M = m + 2$, then $\Gamma_M = X$ and the formula
 (\ref{equChromatic}) is valid by what we just proved.
 
   If $M > m + 2$, we note that $Y = \Gamma_{M - 1}$ and
that the intersection of graphs $X$ and $Y$ is the complete graph $K_{m + 1}$. Then by formula (\ref{equRead}), we have
\bal{
\chi_{\Gamma_M} &= \frac{\chi_X (\lambda) \chi_{\Gamma_{M - 1}}(\lambda)}{\lambda^{(m + 1)}} 
\\
&=\frac{\lambda \cdots (\lambda - m + 1) (\lambda - m)^2  \chi_{\Gamma_{M - 1}}(\lambda)}{\lambda \cdots (\lambda - m)}
\\
&= (\lambda - m) \chi_{\Gamma_{M - 1}}(\lambda),
}
which proves formula (\ref{equChromatic}) by induction supposition. 
\end{proof}

\begin{proof}[Proof of Theorem \ref{theoStairEntropy}]
For large $M$, we have from Theorem \ref{theoStair}, 
\bal{
\frac{\log_2(|\mathcal{T}(St_M^{(n)})|)}{M} &= \frac{\log_2\Gamma\big(\frac{n + 1}{2}\big)}{M} + \frac{(M - (n+1)/2)   \log_2\big(\frac{n + 1}{2}\big)}{M} 
\\
&\longrightarrow \log_2 (n + 1) - 1 , \text{ as } M \to \infty.
}
\end{proof}


%

\appendix
\section{Stanley's results about the number of ribbon tilings}

 \begin{figure}[htbp]
\begin{minipage}[b]{0.42\linewidth}
\centering
              \includegraphics[width=\linewidth]{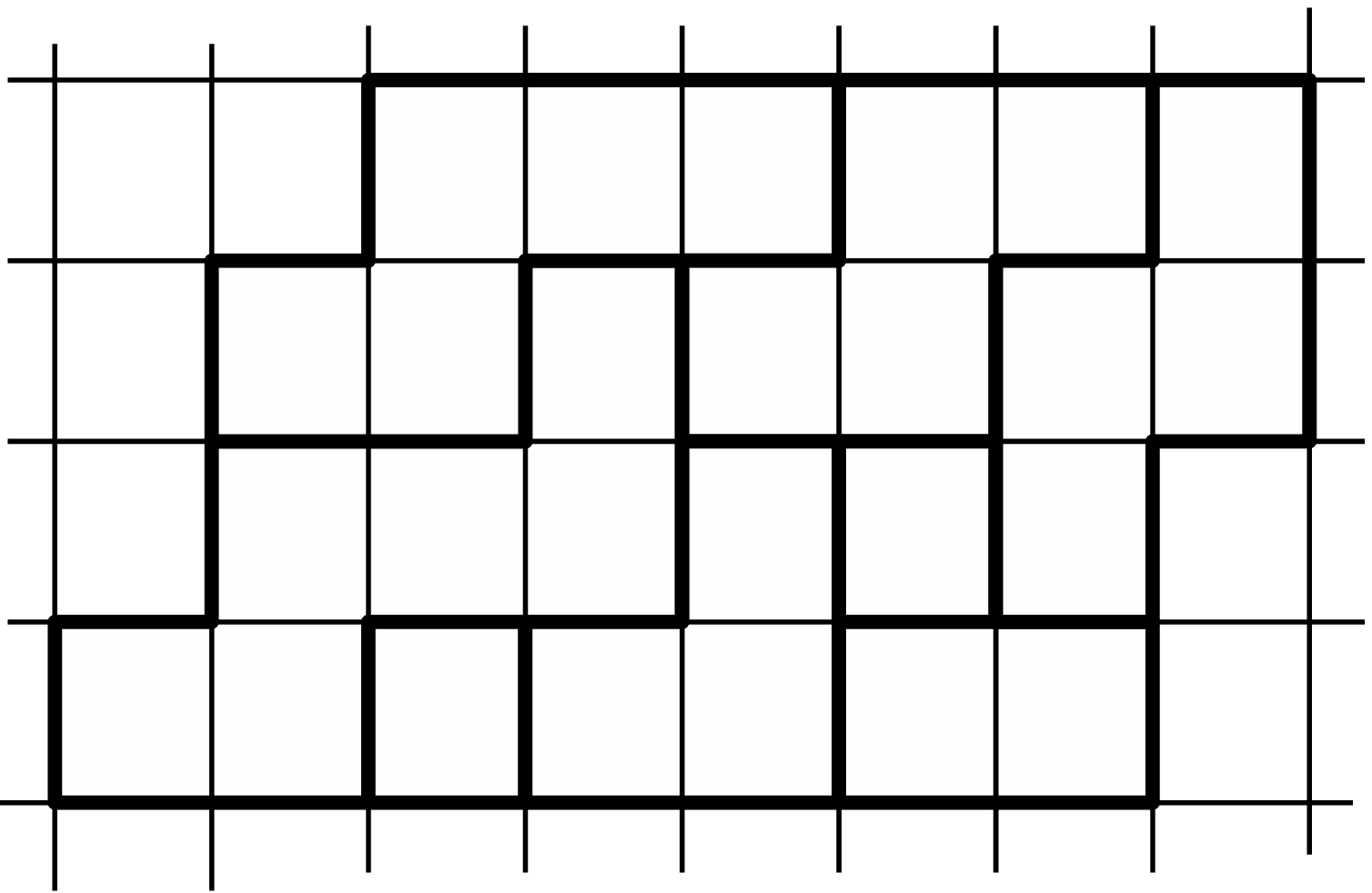}
              \caption{A  ribbon tiling of $8877/211$}
              \label{figStanleyTiling}
\end{minipage}
\hspace{0.5cm}
\begin{minipage}[b]{0.48\linewidth}
\centering
              \includegraphics[width=\linewidth]{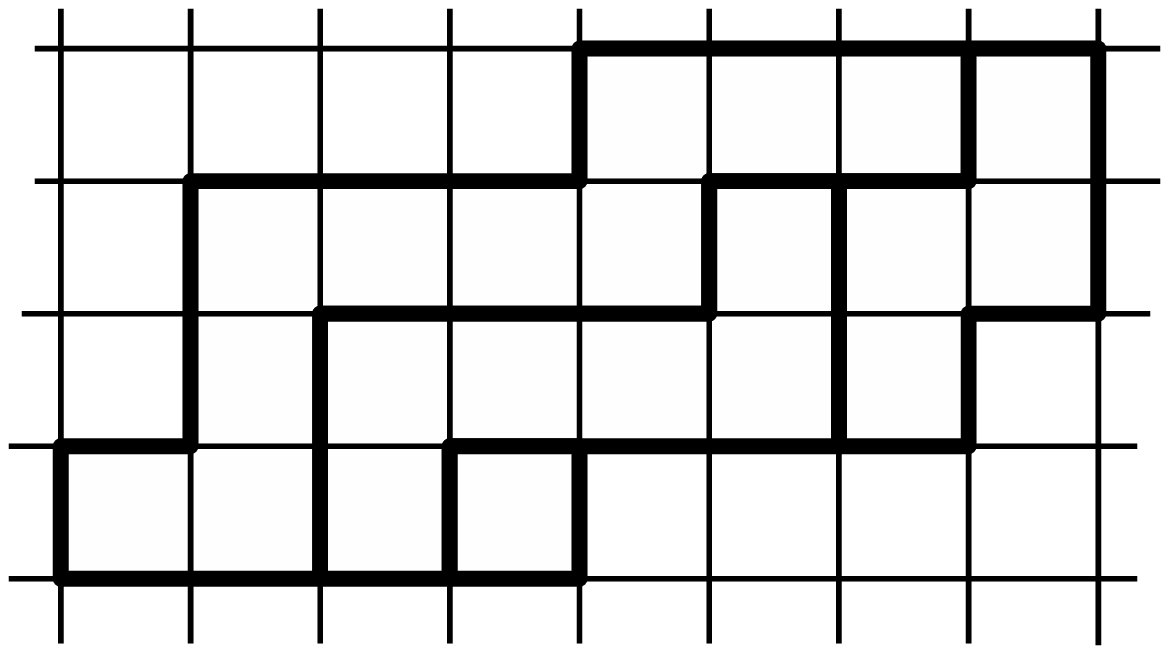}
              \caption{A minimal ribbon tiling of $8874/411$}
              \label{figMinimalRibTiling}
\end{minipage}
\end{figure}

Stanley considers the ribbon tilings of connected skew shapes. (He uses the name \emph{border strip decompositions} for these tilings, see an example in Figure \ref{figStanleyTiling}.)  A skew shape $\lambda/\mu$ is the difference of two Young diagrams $\lambda$ and $\mu$, such that $\mu$ is inside of $\lambda$. In particular, skew shapes include all Young diagrams, and in particular, all rectangles. Stanley works under assumption that lengths of ribbons in such a tiling can be arbitrary, which is different from our assumption of fixed length. 

 In Exercise 7.66 (p. 470, with a solution on p. 521)  in \cite{stanley99}, it is shown that the number of ribbon tilings can be written as a product of certain Fibonacci numbers. 

Here is a summary of Stanley's result applied to an $M$-by-$N$ rectangle with  $N \geq M$. The number of tilings of the rectangle is given by the product:
\bal{
f_{M,N} =\bigg( \prod_{k = 1}^{M - 1}F_{2k+2}^2 \bigg) \Big(F_{2M + 1}\Big)^{N - M},
}
where $F_k$ are Fibonacci numbers, $F_1 = 1, F_2 = 1, F_3 = 2, F_4 = 3, F_5 = 5, \ldots$. (For example, for the $1 \times 2$ rectangle and  the $2 \times 2$ square we have $f_{1,2} = F_3 = 2$ and $f_{2, 2} = F_4^2 = 9$, respectively.)

In particular, for squares, $N = M$, we can write the entropy per unit area as 
\bal{
N^{-2} \log_2(f_{N, N}) = 2 N^{-2} \sum_{k=1}^{N - 1}\log_2(F_{2k + 2}).
}
By using the asymptotic approximation for Fibonacci numbers, the asymptotic expression for the entropy per unit area is
\bal{
N^{-2} \log_2(f_{N, N}) \sim 2 N^{-2}  \log_2 (\phi)\sum_{k=1}^{N - 1}(2k + 2) \sim 2 \log_2(\phi) \approx 1.3885\ldots,
}
where $\phi = (1 + \sqrt{5})/2$ is the golden ratio. It is somewhat difficult to compare this result with our findings, since it not clear what is the average entropy per tile. The number of tiles $T$ is different in different Stanley's tilings of an $N\times N$ rectangle and we do not know how the expectation of $T$ in a random ribbon tiling depends on $N$.

The second result was obtained in \cite{stanley2002}. It is still assumed that the lengths of ribbons in a ribbon tiling are arbitrary. A ribbon tiling is called \emph{minimal} if there is no other tiling with a smaller number of ribbons (see an example in Figure \ref{figMinimalRibTiling}). For skew shapes, Stanley determined the number of tiles in a minimal tiling and the number of minimal ribbon tilings. If we specialize his results to $M$-by-$N$ rectangles with $M \leq N$, then the number of ribbons in a minimal tiling always equals $M$. 
(So, for large rectangles, most of the tiles in a minimal tiling are long, with the average length equal $N$.) In the case of $M\times N$ rectangles with $M \leq N$,  Stanley's formula for the number of minimal tilings reduces to
\bal{
h_{M,N} = (M!)^2.
}
In particular,  the  per-tile entropy is 
\bal{
M^{-1} \log_2(h_{M, N}) = 2 M^{-1} \log_2 M! \sim 2 \log_2 (M/e).
}

As a consequence, if $M$ is fixed and $N$ is growing then the asymptotic per-tile entropy does not depend on the average length of the tile $N$. This is in contrast to our results, where the per-tile entropy grows at the logarithmic rate in the length of the tile. 

If $M = N$ and both are growing, then the per-tile entropy grows at the logarithmic rate in the average length of the tile $N$, similar to our results.

\subsection*{Acknowledgements}
The second author was supported by Simons Foundation through the program ``Mathematics and Physical Sciences--Travel Support for Mathematicians'', Award ID 523587. 
%
%
%
%


\end{document}